\documentclass[12pt]{article}

\usepackage{amsmath,amssymb}

\begin{document}

\def \1{{\bf 1}}
\def \a{{{\frak a}}}
\def \al{\alpha}
\def \ar{{\alpha_r}}
\def \A{{\mathbb A}}
\def \Ad{{\rm Ad}}
\def \b{{{\frak b}}}
\def \bs{\backslash}
\def \B{{\cal B}}
\def \cent{{\rm cent}}
\def \C{{\mathbb C}}
\def \CA{{\cal A}}
\def \CB{{\cal B}}
\def \CE{{\cal E}}
\def \CF{{\cal F}}
\def \CG{{\cal G}}
\def \CH{{\cal H}}
\def \CL{{\cal L}}
\def \CM{{\cal M}}
\def \CN{{\cal N}}
\def \CP{{\cal P}}
\def \CQ{{\cal Q}}
\def \CO{{\cal O}}
\def \CS{{\cal S}}
\def \det{{\rm det}}
\def \diag{{\rm diag}}
\def \dom{{\rm dom}\hspace{2pt}}
\def \e{\epsilon}
\def \End{{\rm End}}
\def \Fx{{\frak x}}
\def \FX{{\frak X}}
\def \g{{{\frak g}}}
\def \ga{\gamma}
\def \Ga{\Gamma}
\def \h{{{\frak h}}}
\def \Hom{{\rm Hom}}
\def \im{{\rm im}\hspace{2pt}}
\def \Ind{{\rm Ind}}
\def \k{{{\frak k}}}
\def \K{{\cal K}}
\def \la{\lambda}
\def \lap{\triangle}
\def \Lie{{\rm Lie}\hspace{2pt}}
\def \m{{{\frak m}}}
\def \mathqed{\hfill \Box} 

\def \mod{{\rm mod}}
\def \n{{{\frak n}}}
\def \name{\bf}
\def \N{\mathbb N}
\def \ord{{\rm ord}}
\def \O{{\cal O}}
\def \p{{{\frak p}}}
\def \ph{\varphi}
\def \prf{{\bf Proof: }}
\def \qed{$ $\newline $\frac{}{}$\hfill {\rm Q.E.D.}\vspace{15pt}}
\def \Q{\mathbb Q}
\def \res{{\rm res}}
\def \R{{\mathbb R}}
\def \Re{{\rm Re \hspace{1pt}}}
\def \ra{\rightarrow}
\def \rank{{\rm rank}}
\def \Rep{{\rm Rep}}
\def \sign{{\rm sign}\hspace{2pt}}
\def \supp{{\rm supp}}
\def \t{{{\frak t}}}
\def \T{{\mathbb T}}
\def \tr{{\hspace{1pt}\rm tr\hspace{1pt}}}
\def \vol{{\rm vol}}
\def \V{{\cal V}}
\def \z{\frak z}
\def \Z{\mathbb Z}
\def \={\ =\ }

\newcommand{\rez}[1]{\frac{1}{#1}}
\newcommand{\der}[1]{\frac{\partial}{\partial #1}}
\newcommand{\norm}[1]{\parallel #1 \parallel}
\newcommand{\krn}[3]{\langle #1 | #2 | #3 \rangle}

\newcounter{lemma}
\newcounter{corollary}
\newcounter{proposition}
\newcounter{theorem}
\newcounter{zwisch}

\renewcommand{\subsection}{\refstepcounter{subsection}\stepcounter{lem 
ma} 

	\stepcounter{corollary} \stepcounter{proposition}
	\stepcounter{conjecture}\stepcounter{theorem}
	$ $ \newline
	{\bf \arabic{section}.\arabic{subsection}\hspace{8pt}}}

\newtheorem{conjecture}{\stepcounter{lemma} \stepcounter{corollary} 	
	\stepcounter{proposition}\stepcounter{theorem}
	\stepcounter{subsection}\hskip-12pt Conjecture}[section]
\newtheorem{lemma}{\stepcounter{conjecture}\stepcounter{corollary}	
	\stepcounter{proposition}\stepcounter{theorem}
	\stepcounter{subsection}\hskip-7pt Lemma}[section]
\newtheorem{corollary}{\stepcounter{conjecture}\stepcounter{lemma}
	\stepcounter{proposition}\stepcounter{theorem}
	\stepcounter{subsection}\hskip-7pt Corollary}[section]
\newtheorem{proposition}{\stepcounter{conjecture}\stepcounter{lemma}
	\stepcounter{corollary}\stepcounter{theorem}
	\stepcounter{subsection}\hskip-7pt Proposition}[section]
\newtheorem{theorem}{\stepcounter{conjecture} \stepcounter{lemma}
	\stepcounter{corollary}\stepcounter{proposition}		
	\stepcounter{subsection}\hskip-11pt Theorem}[section]

\title{Harmonic Analysis on the quotient  
$\Q^\times\ltimes\Q\bs\A^1\ltimes \A$}
\author{{\small by}\\ {} \\ Anton Deitmar}
\date{}
\maketitle

\pagestyle{myheadings}
\markright{HARMONIC ANALYSIS ON THE QUOTIENT...}

\noindent
{\small {\bf Abstract. }
Let $G$ be the semidirect product $\A^1\ltimes \A$ of the adeles and  
the norm $1$ ideles.
Let $\Ga$ be the discrete subgroup $\Q^\times\ltimes\Q$.
In this paper the trace formula for this setting is established and  
used to give the complete decomposition of the $G$-representation on  
$L^2(\Ga\bs G)$.
It turns out that every character of the norm-$1$ idele class group  
gives a one dimensional isotype and the complement of those consists  
of one irreducible representation.}

\tableofcontents

$$ $$

\noindent
{\bf Introduction. }
The group of norm-$1$ ideles $\A^1$ over $\Q$ acts by multiplication  
on the additive group of adeles $\A$.
Let $G=\A^1\ltimes\A$ be their semidirect product.
Let $\Ga$ be its arithmetic subgroup $\Q^\times\ltimes\Q$.
In \cite{gold} D. Goldfeld gave a Schwartz kernel on $\Ga\bs G$ whose  
geometrical trace equals the prime number side of the explicit  
formulas of A. Weil.
A spectral interpretation in the sense of the Selberg trace formula  
is called for.

In this paper we give a classification of the irreducible unitary  
representations of $G$ and describe their trace functionals as  
integrals on the group.
Next we set up the trace formula and make the geometric side explicit  
by evaluating the corresponding orbital integrals.
Comparing these with the traces yields a full description of the  
``automorphic spectrum'' of $G$, i.e. a complete decomposition of the  
$G$-representation on $L^2(\Ga\bs G)$.
It turns out that every character of the norm-$1$ idele class group  
$\A^1 /\Q^\times$ gives a one dimensional isotype.
The joint orthogonal complement of those forms a single irreducible  
representation of $G$.

This case is intermediate between the (from the point of  
representation theory) simple case of $GL_1(\A)$ and the complicated  
case $GL_2(\A)$, where no explicit description of the automorphic  
spectrum is available.

$$ $$

\section{The group $G=\A^1\ltimes\A$ and its Hecke algebra}
Let  $\A$ the ring of adeles over $\Q$ and let $\A^\times$ the group  
of ideles.
Let $\A_{fin}$ be the ring of finite adeles, then $\A  
=\A_{fin}\times\R$.
At any rational prime $p$ let the absolute value $|.|_p :  
\Q_p\ra\R_{\ge 0}$ be normalized so that $|p|_p=p^{-1}$.
The product formula tells us that $\Q^\times$ embedded diagonally  
into $\A^\times$ lies in the set $\A^1$ of ideles $a$ with $|a|=1$.
The multiplicative group $\A^1$ acts on the additive group $\A$ by  
multiplication and we can form the semidirect product
$$
G\= \A^1\ltimes \A.
$$
We will write the product of $(a,x)$ and $(b,y)$ in $G$ as
$$
(a,x)(b,y) \= (ab ,bx+y).
$$
The first simple observation on the structure of the group $G$ is

\begin{lemma}
The center of $G$ is trivial.
\end{lemma}

On $\A$ we fix a Haar measure $dx=\otimes_vdx_v$ where $dx_v$ is  
normalized by $\vol(\Z_p)=1$ if $v=p$ and $dx_\infty$ equals the  
Lebesgue measure on $\R$.
We fix local multiplicative Haar measures $d^\times a_v$ by  
$\vol(\Z_p^\times)=1$ if $v=p$ and  
$d^*a_\infty=\frac{da_\infty}{|a_\infty|_\infty}$.
Then $d^\times a =\prod_v d^*a_v$ forms a Haar measure of  
$\A^\times$.

Let $i:\A^1\hookrightarrow \A^\times$ be the inclusion and let  
$s:\R_+^\times \ra \A^\times$ given by the embedding of $\R$ at  
infinity then the map $(i,s):\A^1\times\R_+^\times\ra\A^\times$ is an  
isomorphism.
The Haar measure $d^\times a$ thus factors into $d^\times a = d^1  
a\otimes\frac{dt}{t}$ and fixes a Haar measure $d^1a$ on $\A^1$.

Next let $x\in\A$ and let $\A^1_x$ be its stabilizer group in $\A^1$,  
i.e. the set of all $a\in\A^1$ such that $ax=x$.
Then $\A^1_x$ is the intersection of $\A^1$ with $\A^\times_x$, which  
is the restricted product of all completions $\Q_v$ such that  
$x_v=0$.
Taking the product over all $d^\times a_v$ gives a canonical Haar  
measure on $\A^\times_x$.
Let $N_x$ be the image of $\A^\times_x$ under $|.|$ then we get a  
product decomposition $\A^\times_x=\A^1_x\times N_x$. 

Now $N_x$ either is trivial or equals $\R^\times_{>0}$ or is  
countable.
In either case there is a natural choice for a Haar measure on $N_x$  
which then gives a natural choice for a Haar measure on $\A^1_x$.

\begin{lemma}
The group $G$ is unimodular.
The measure $d^1a\otimes dx$ is a Haar measure on $G$.
\end{lemma}

\prf
The groups $\A^1$ and $\A$ are unimodular and $\A^1$ acts measure  
preservingly on $\A$, hence $\A^1\ltimes \A$ is unimodular  
\cite{fell-doran}.
The second fact is easily proven.
\qed

Let $\CH =\CH(G)$ be the {\bf Hecke algebra} of $G$, i.e. the  
convolution algebra of all $\C$-valued compactly supported smooth  
functions on $G=\A^1\ltimes\A$.
Smooth here means infinitely often differentiable along the real  
coordinate and locally constant in the $p$-adic directions.
We will make this more precise.
Let $\hat{\Z} =\prod_p\Z_p\subset\A$ be the profinite completion of  
$\Z$ and $\hat{\Z}^\times =\prod_p\Z_p^\times\subset\A^1$ its unit  
group.
Then $\A^1 =\Q^\times\hat{\Z}^\times$ and the product is direct.
This gives an isomorphism
$$
G\ \cong\ \Q^\times \ltimes\left( \hat{\Z}^\times\ltimes\A\right).
$$
A function $f$ on $\hat{\Z}^\times\ltimes\A$ is in the Hecke algebra  
of $\hat{\Z}^\times\ltimes\A$ if $f$ is a finite sum of functions of  
the form $f=(\prod_p f_p)f_\infty$, where $f_\infty\in  
C_c^\infty(\R)$ and $f_p\in C_c^\infty(\Z_p^\times\ltimes\Q_p)$ such  
that $f_p$ is the characteristic function of ${\Z_p^\times\ltimes  
\Z_p}$ for almost all $p$.
Now a function $h$ on $G=\Q^\times \ltimes\left(  
\hat{\Z}^\times\ltimes\A\right)$ is in the Hecke algebra $\CH  
=\CH(G)$ if $h$ is a finite sum of functions of the form $h=\1_{\{  
r\} }f$, where $r$ is in $\Q^\times$ and $f$ is in the Hecke algebra  
of $\hat{\Z}^\times\ltimes\A$.

\section{The unitary dual}
We fix a canonical basic character of $\A$ as   
$\psi=(\prod_p\psi_p)\psi_\infty$ with
$\psi_p(\Z_p)=1$, $\psi_p(p^{-n})=e^{2\pi i/p^n},$
and $\psi_\infty(x)=e^{2\pi ix}.$
Note that $\psi$ is chosen in a way that the lattice $\Q\subset\A$ is  
its own dual, i.e. 

$$
\langle x,y\rangle\=1\ \ \forall y\in \Q\hspace{20pt}  
\Leftrightarrow\hspace{20pt} x\in \Q.
$$
For any $y\in\A$ we write $\psi_y$ for the character  
$\psi_y(x)=\psi(xy)$.
The group $\A^1$ acts on $\hat{\A}$ via $a\psi_x =\psi_{a^{-1}x}$.

Let $x\in\A$ then its stabilizer group $\A^1_x =\{ a\in\A^1 | ax=x\}$  
is the group of norm 1 elements in the restricted product of the  
$k_v$ over all places $v$ such that $x_v=0$.
Let $S$ denote this set of places, let $\A_S^\times$ be the  
restricted product over all places in $S$ and let $\A^{\times,S}$ be  
the restricted product over all places outside $S$.
Then $\A^\times =\A_S^\times\times\A^{\times ,S}$ and $ \A^1 = \{  
x=(x_S,x^S) : |x_S||x^S|=1\}$.
Let $V = |\A^\times_S|\cap |\A^{\times,S}|\ \subset\ \R_+^\times$ the  
joint value group, then $\A^1 = \A^1_S \times \A^{1,S}\times V$, so  
$\A_S^1 =\A^1_x$ is a direct summand of $\A^1$.
It follows that any character $\alpha$ of $\A^1_x$ can be lifted to a  
character of $\A^1$ which will be denoted by the same letter.
Fix a Haar-measure $\mu_x$ on $\A^1/\A^1_x$ and
let $H_x$ be the Hilbert space $L^2(\A^1/\A^1_x)$.
We define a representation $\pi_x$ of  $G=\A^1\ltimes \A$ on $H_x$ by
$$
\pi_x((a,y)) \ph(b) \= \psi((ab)^{-1} xy)\ph(ab),\ \ \ \ph\in H_x.
$$
This gives a unitary representation $(\pi_x,H_x)$ of $G$.
Let $\alpha$ be a unitary character of $\A^1_x$ lifted to $G$ then  
$\pi_x\otimes\alpha$ again is a unitary representation of $G$.
This defines a map 

$$
F:\bigcup_{[x]\in\A /\A^1}\hat{\A^1_x} \ra \hat{G}.
$$

\begin{theorem}
The map $F$ is a bijection, i.e. the unitary dual of $G$ is  
parametrized by $\bigcup_{[x]\in\A /\A^1}\hat{\A^1_x}$.
\end{theorem}

\prf
Take any $\pi\in\hat{G}$, then the restriction $\pi|_{\A}$ can be  
disintegrated so that there is a measurable bundle $\tau\mapsto  
H_\tau$ of Hilbert spaces over $\hat{\A}$ and a measure $\mu$ on  
$\hat{\A}$ such that
$$
\pi|_\A\= \int_{\hat{\A}} H_\tau d\mu(\tau).
$$
The group $\A^1$ acts on this space in a way that $a\in\A^1$ maps  
$H_\tau$ to $H_{^a\tau}$.
Thus irreducibility of $\pi$ implies that $\mu$ is supported on a  
single $\A^1$-orbit in $\hat{\A}$.
Now fix $\tau =\psi_x$, $x\in\A$ in that orbit, then the stabilizer  
$\A^1_x$ acts on $H_{\psi_x}$.
Let $V\subset H_{\psi_x}$ be a proper $\A^1_x$-stable  subspace then
$\int_{\hat{\A}}V_\tau d\mu(\tau)$ is a proper $G$-subrepresentation,  
where $V_{\psi_{ax}}=aV$.
Thus irreducibility forces $H_\tau$ to be irreducible, that is, one  
dimensional and $\A^1_x$ acts on it by a character $\alpha$, the  
surjectivity of $F$ follows.
Backwards the argument also shows that $\pi\in\hat{G}$ uniquely gives  
the orbit $[x]$ and the character $\alpha$.
The theorem follows.
\qed

For $\pi\in\hat{G}$ we will now give a formula for the trace  
distribution $h\mapsto\tr \pi(h)$.
Let $x\in\A$ and $\pi =\pi_x$.
Let $\alpha$ be a unitary character of $\A_x^1$ lifted to $G$.
We will consider the representation $\pi_x\otimes\alpha$.
The representation space is $H_x=L^2(\A^1/\A^1_x)$. For $\ph\in H_x$  
and $h\in\CH$ we have that $\pi_x\otimes\alpha(h)\ph (a)$ equals
\begin{eqnarray*}
 &{}& \int_G h(z)\pi_x(z)\alpha(z)\ph(a) dz\\
	&=& \int_{\A^1}\int_\A  
h(b,y)\psi((ab)^{-1}xy)\alpha(b)\ph(ab)dyd^\times b\\
	&=& \int_{\A^1/A^1_x}\int_{A^1_x}\int_\A  
h(bb',y)\psi((ab)^{-1}xy)\alpha(bb')\ph(ab) dy d^\times b' d^\times  
b\\
	&=& \int_{\A^1/A^1_x}\int_{A^1_x}\int_\A  
h(a^{-1}bb',y)\psi(b^{-1}xy)\alpha(a^{-1}bb')\ph(b) dy d^\times b'  
d^\times b
\end{eqnarray*}
so that $\pi_x\otimes\alpha(h)$ has kernel
$$
k(a,b)\= \int_{\A^1_x}\int_\A  
h(a^{-1}bb',y)\alpha(a^{-1}bb')\psi(b^{-1}xy)dyd^\times b',
$$
which is smooth and compactly supported so that
\begin{eqnarray*}
\tr\pi_x\otimes\alpha(h) &=& \int_{\A^1 /\A^1_x} k(a,a)d^\times a\\
	&=& \int_{\A^1 /\A^1_x} \int_{\A^1_x}\int_\A  
h(b',y)\alpha(b')\psi(a^{-1}xy) dy d^\times b'd^\times a.
\end{eqnarray*}
Let $\hat{h}(a,x) = \int_\A h(a,y)\psi(xy)dy$ the $\A$-Fourier  
transform.
For two functions $f,g$ on $\A^1_x$ let
$$
\langle f,g\rangle_{\A^1_x} \= \int_{\A^1_x} f(a)g(a)d^\times a.
$$
We have shown

\begin{proposition} \label{trace}
For $h\in\CH$ and $x\in\A$ we have
$$
\tr\pi_x\otimes\alpha(h) \= \int_{\A^1 /\A^1_x}  
\langle\alpha,\hat{h}(.,ax)\rangle_{\A_x^1} d^\times a.
$$
\end{proposition}

We will consider some special cases.
At first let $x=1$ and note

\begin{lemma}\label{trpi1}
For $h\in\CH$ we have
\begin{eqnarray*}
\tr\pi_1(h) &=& \int_{\A^1} \hat{h}(1,a) d^\times a\\
	&=& \sum_{q\in\Q^\times} \int_{\hat{\Z}^\times}  
\hat{h}(1,qa)d^\times a.
\end{eqnarray*}
\end{lemma}

\prf
This follows from the proposition and $\A^1 =\Q^\times  
\hat{\Z}^\times$.
\qed

Next write $K$ for the compact group  
$K=\A^1/\Q^\times\cong\hat{\Z}^\times$.
Note that the natural map $G=\A^1\ltimes\A\ra\A^1\ra\A^1/\Q^\times  
=K$ gives $\hat{K}\hookrightarrow \hat{G}$.

Let $R_K$ denote the right regular representation of $K$ or $G$ on  
$L^2(K)$.

\begin{lemma}\label{trrk}
For any $h\in\CH$ we have
$$
\tr R_K(h) \= \sum_{q\in\Q^\times} \hat{h}(q,0).
$$
\end{lemma}

\prf
Let $h_K(a) = \sum_{q\in\Q^\times}\int_\A h(qa,x) dx$.
Then $R_K(h)=R_K(h_K)$ where on the right hand side $R_K$ is  
considered as a representation of $K$.
This operator has a smooth kernel on $K$ given by
$$
k(a,b) \= h_K(a^{-1}b).
$$
Therefore its trace equals
$$
\int_K k(a,a) d^\times a\= h_K(1)\= \sum_{q\in\Q^\times}\hat{h}(q,0).
$$
\qed

\section{The decomposition of $L^2(\Ga\bs G)$}
The group $\Ga =\Q^\times\ltimes \Q$ forms a cocompact lattice in  
$G$.
We consider the unitary representation $R$ of $G$ on $L^2(\Ga\bs G)$.
So for $\ph\in L^2(\Ga\bs G)$ we have $R(y)\ph(x)=\ph(xy)$.

Here comes the main result of this paper.

\begin{theorem}
As a $G$-representation we have
\begin{eqnarray*}
L^2(\Ga\bs G) &\cong& R_K \oplus\pi_1\\
	&\cong& \bigoplus_{\chi\in\hat{K}}\chi \oplus\pi_1
\end{eqnarray*}
\end{theorem}

\prf
At first we need to know the following

\begin{lemma}
The space $L^2(\Ga\bs G)$ splits into a discrete sum of irreducible  
representations with finite multiplicities. More precisely
$$
L^2(\Ga\bs G) \= \bigoplus_{\pi\in\hat{G}}N(\pi)\pi,
$$
where $N(\pi)\in\N_0$ and $N(\pi)=0$ for all but countably many  
$\pi$.
\end{lemma}

\prf
For any $h\in\CH$ the operator $R(h) =\int_Gh(x)R(x)dx$ is a bounded  
linear operator on $L^2(\Ga\bs G)$.
On the Hecke algebra $\CH$ we have the involution  
$h^*(x)=\overline{h(x^{-1})}$.
Assume $h=h^*$ then $R(h)$ is selfadjoint.
For $\ph\in L^2(\Ga\bs G)$ we have
\begin{eqnarray*}
R(h)\ph(x) &=& \int_G h(y)\ph(xy) dy\\
	&=& \int_G h(x^{-1}y)\ph(y) dy\\
	&=& \sum_{\ga\in\Ga}\int_\CF h(x^{-1}\ga y)\ph(y) dy\\
	&=& \int_{\Ga\bs G}k_h(x,y)\ph(y) dy,
\end{eqnarray*}
where $\CF$ is an arbitrary fundamental domain for $\Ga$ in $G$ and  
$k_h(x,y)=\sum_{\ga\in\Ga}h(x^{-1}\ga y)$.
The sum is locally finite so $k_h$ is continuous on the compact space  
$\Ga\bs G\times \Ga\bs G$.
This implies that $R(h)$ is a Hilbert-Schmidt operator, hence  
compact, so has discrete spectrum with finite multiplicities away  
from zero.
The Hecke algebra $\CH$ contains an approximate identity of  
selfadjoint elements so the lemma follows.
\qed

\begin{lemma}
For any $h\in\CH$ the operator $R(h)$ is of trace class. Its trace  
equals
$$
\sum_{\pi\in\hat{G}} N(\pi)\tr\pi(h)\= h(1)  
+\int_{\A^1}h(1,a)d^\times a +\sum_{\alpha\in\Q^\times -\{ 1\}  
}\int_\A h(\alpha ,x) dx.
$$
\end{lemma}

\prf
This is the Selberg trace formula for the group $G$.
The only point that goes beyond standard considerations is the  
determination of the orbital integrals.

Let $h\in \CH$.
The sum defining $k_h$ in the last proof being locally finite implies  
$k_h$ is smooth.
Hence the operator $R(h)$ is of trace class, its trace being
$$
\tr R(h) \= \sum_{\pi\in\hat{G}} N(\pi) \tr\pi(h).
$$
On the other hand the trace equals the integral over the diagonal of  
the kernel, so
\begin{eqnarray*}
\tr R(h) &=& \int_{\Ga\bs G}k_h(x,x) dx\\
	&=& \int_{\Ga \bs G}\sum_{\ga\in\Ga}h(x^{-1}\ga x)dx\\
	&=& \sum_{\ga\in\Ga}\int_{\CF}h(x^{-1}\ga x)dx,
\end{eqnarray*}
where $\CF$ is a fundamental domain for the $\Ga$-action on $G$.
Next we split the sum over $\Ga$ into a sum over all conjugacy  
classes $[\ga]$.
Let $\Ga_\ga$ be the centralizer of $\ga$ in $\Ga$ and $G_\ga$ the  
centralizer in $G$. We obtain
\begin{eqnarray*}
&{}& \sum_{[\ga]}\sum_{\sigma\in\Ga_\ga\bs\Ga} \int_\CF h((\sigma  
x)^{-1} \ga \sigma x)dx\\
&=& \sum_{[\ga]} \int_{\Ga_\ga\bs G} h(x^{-1} \ga x)dx\\
&=& \sum_{[\ga]} \vol(\Ga_\ga\bs G_\ga)\int_{G_\ga\bs G} h(x^{-1} \ga  
x)dx
\end{eqnarray*}
for any choice of Haar measures on the centralizers $G_\ga$.
Suppose such measures chosen let
$\CO_\ga(h) \= \int_{G_\ga\bs G} h(x^{-1} \ga x)dx$ be the {\bf  
orbital integral}.

We now shall determine the conjugacy classes in $\Ga$ and their  
orbital integrals.
We start with $\ga =e$, the trivial element.
Then $G_\ga =G$ and $\vol(\Ga_\ga\bs G_\ga)=\vol(\Ga\bs G) =1$ with  
the Haar measure already chosen.
Further $\CO_\ga (h) = h(1)$.
Next let $\ga =(a_\ga ,x_\ga)\in\Ga$ and $z=(a,x)\in G$ then
$$
z^{-1}\ga z \= (a_\ga, x(1-a_\ga) +ax_\ga).
$$
First consider the case $a_\ga =1$.
Then $z^{-1}\ga z=(1,ax_\ga)$.
The case $x_\ga=0$ has already been dealt with.
So let $x_\ga\ne 0$.
Then the centralizer $G_\ga$ equals the set of all $(1,x)$, $x\in\A$  
and so $G_\ga\cong\A$.
We choose the Haar measure $dx$ on $G_\ga$.
We see that all $(1,x_\ga)$ with $x_\ga\in\Q^\times$ form a single  
conjugacy class in $\Ga$.
This gives the second summand in the theorem.

Now assume $a_\ga\ne 1$.
We get
$$
(1,\frac{x_\ga}{a_\ga-1})^{-1}  
(a_\ga,x_\ga)(1,\frac{x_\ga}{a_\ga-1})\= (a_\ga ,0),
$$
which implies that the set of $(a_\ga,x_\ga)$ for a fixed $a_\ga\ne  
1$ form a single conjugacy class. We choose the representative $\ga  
=(a_\ga ,0)$, then we get
$G_\ga\cong\A$ and
$$
\int_{G_\ga\bs G}h(z^{-1}\ga z)dz\= \int_\A h(a_\ga ,x(1-a_\ga))dx,
$$
which gives the summand with $\alpha =a_\ga$.
The lemma follows.
\qed

With the aid of this we prove

\begin{lemma}
For $h\in\CH$ we have
$$
\tr R(h)\= \tr R_K(h) +\tr\pi_1(h).
$$
\end{lemma}

\prf
The Lemmas \ref{trpi1} and \ref{trrk} give that the right hand side  
equals
$$
\sum_{q\in\Q^\times}\hat{h}(q,0) + \int_{\A^1}\hat{h}(1,a)d^\times a.
$$
The second summand equals
\begin{eqnarray*}
\sum_{q\in\Q^\times}\int_{\hat{\Z}^\times} \hat{h}(1,qa)d^\times a	 
&=& \int_{\hat{\Z}^\times}\sum_{q\in\Q^\times} \hat{h}(1,qa)d^\times  
a\\
	&=& \int_{\hat{\Z}^\times}\sum_{q\in\Q^\times}  
h(1,qa)-\hat{h}(1,0)+h(1,0)d^\times a
\end{eqnarray*}
according to the Poisson Summation formula for $\A$ and the lattice  
$\Q$.
This implies that $\tr R_K(h)+\tr\pi_1(h)$ equals
$$
\sum_{q\in\Q^\times-\{ 1\} }\hat{h}(q,0) +h(1)  
+\int_{\A^1}h(1,a)d^\times a.
$$
\qed

A unitary representation $\pi$ is called {\bf traceable} if $\pi(h)$  
is of trace class for every $h\in\CH$.
It follows that $\pi$ is a countable sum of irreducibles.
To finish the prove of the theorem it remains to show

\begin{lemma}
Let $R,S$ be two traceable representations of $G$ and assume that  
$\tr R(h)=\tr S(h)$ for every $h\in\CH$. Then $R\cong S$.
\end{lemma}

\prf
An element $x$ of $\A$ is called {\bf regular} if for every place $v$  
the entry $x_v$ is nonzero. This is equivalent to $\A^\times_x$ being  
trivial.
Otherwise $x$ is called {\bf irregular}.
The representation $\pi_x$ is called regular if $x$ is. Let  
$\hat{G}_{reg}$ be the set of regular representations.
Write the decomposition of $R$ into irreducibles as  
$R=\bigoplus_{\pi\in\hat{G}} N_R(\pi)\pi$ and let  
$R_{reg}=\bigoplus_{\pi\in\hat{G}_{reg}} N_R(\pi)\pi$ be the regular  
part.
Define $S_{reg}$ analogously.
Consider $h\in\CH$ of the form $h=h_1\otimes h_2$ with  
$h_1\in\CH(\A^1)$ and $h_2\in\CH(\A)$.
Assume $h_1(1)=0$, then $\tr\pi(h)=0$ for any regular $\pi$ as  
Proposition \ref{trace} shows.
There still is enough freedom in the choice of $h$ to conclude
$\tr R_{reg}(h)=\tr S_{reg}(h)$ for any $h\in\CH$.
Let $R_{reg}$ be given by the sequence $(x_n)$ and $S_{reg}$ by  
$(y_n)$.
Then
$$
\tr R_{reg}(h)\= \sum_n\int_{\A^1}\hat{h}(1,ax_n)d^\times a.
$$
The $x_n$ are determined up to $\A^1$-multiplication only.
The $\A^1$-orbits of the $x_n$ cannot accumulate since then $R(h)$  
would not be of trace class for every $h\in\CH$.
But this means they are discrete and to any given $n$ there is  
$h\in\CH$ such that $\tr\pi_{x_n}(h)\ne 0$ whereas  
$\tr\pi_{x_m}(h)=0$ whenever $x_m\ne x_n$.
Since also the set of $\A^1$-orbits of $(x_j)$ and $(y_j)$ together  
doesn't accumulate it follows that we also can assume  
$\tr\pi_{y_k}(h)=0$ whenever $y_k\ne x_n$.
This leads to 

$$
R_{reg}\ \cong\ S_{reg}.
$$
For the rest we assume that $R$ and $S$ are totally irregular.
Let $\hat{G}_{irr} = \hat{G}-\hat{G}_{reg}$ then $\hat{G}_{irr}$ is  
parametrized by the disjoint union $\bigcup_{x\in(\A  
/\A^1)_{irr}}\hat{\A^1_x}$.
Each $\hat{\A^1_x}$ has a natural topology being a character group.
This defines a topology on the union $\hat{G}_{irr}$.
Again Proposition \ref{trace} shows that the members of $R$ cannot  
accumulate in $\hat{G}_{irr}$. The proof now proceeds as in the  
regular case.

The lemma and the theorem are proven.
\qed

\tiny
\hspace{-20pt}
Math. Inst.\\ INF 288\\ 69120 Heidelberg\\ GERMANY

\end{document}